\documentclass[12pt]{article}
\usepackage{authblk}
\usepackage{amsfonts,amsmath,amssymb,amsthm,mathrsfs,enumerate}

\usepackage{amsmath,bm}
\usepackage{graphicx}
\usepackage{color}
\usepackage{ifpdf}
\usepackage{epstopdf}
\usepackage{subfigure}
\usepackage[numbers,sort&compress]{natbib}
\usepackage{indentfirst}
\usepackage{tipa}
\usepackage{mathtools} 

\usepackage{tabularx}

\newcommand{\Rho}
{\mathbb{P}}

\usepackage{enumitem} 
\setlist[enumerate]{itemsep=0mm}

\usepackage{graphicx,color, booktabs,multirow}
\usepackage{tikz,chemfig}
\usepackage{graphicx,booktabs,multirow,mathrsfs}

\def \vect {\vec{T}}
\def \sett {{\cal T}}

\usepackage{float} 

\usepackage{calc}
\usetikzlibrary{backgrounds}
\usetikzlibrary{shadings,intersections}
\usetikzlibrary{arrows}
\usetikzlibrary{shapes,shapes.geometric,shapes.misc}
\pgfdeclarelayer{edgelayer}
\pgfdeclarelayer{nodelayer}
\usetikzlibrary{automata}
\pgfsetlayers{background,edgelayer,nodelayer,main}
\tikzstyle{none}=[inner sep=0mm]
\tikzstyle{every loop}=[]

\tikzstyle{dotted}=[dash pattern=on \pgflinewidth off 2pt]
\tikzstyle{dashed}=[dash pattern=on 3pt off 3pt]

\usetikzlibrary{backgrounds}
\usetikzlibrary{arrows}
\usetikzlibrary{shapes,shapes.geometric,shapes.misc}
\pgfdeclarelayer{edgelayer}
\pgfdeclarelayer{nodelayer}
\pgfsetlayers{background,edgelayer,nodelayer,main}
\tikzstyle{none}=[inner sep=0mm]
\tikzstyle{every loop}=[]

\tikzstyle{dotted}=[dash pattern=on \pgflinewidth off 2pt]
\tikzstyle{dashed}=[dash pattern=on 3pt off 3pt]

\tikzstyle{new style 0}=[fill=black, draw=black, shape=circle]
\tikzstyle{red style 1}=[fill=red, draw=black, shape=circle]
\tikzstyle{blue style 2}=[fill=blue, draw=black, shape=circle]
\tikzstyle{white style 4}=[fill=white, draw=black, shape=circle]
\tikzstyle{bklack style 5}=[fill=black, draw=black, shape=rectangle]
\tikzstyle{red style 3}=[fill=red, draw=black, shape=rectangle]
\tikzstyle{yellow style 7}=[fill=yellow, draw=black, shape=rectangle]
\tikzstyle{new style 8}=[fill={rgb,255: red,0; green,132; blue,0}, draw={rgb,255: red,0; green,131; blue,0}, shape=circle]

\tikzstyle{new edge style 0}=[-]
\tikzstyle{new edge style 1}=[-, draw=red]
\tikzstyle{new edge style 2}=[-, draw=blue]
\tikzstyle{new edge style 3}=[-, draw={rgb,255: red,0; green,156; blue,0}]

\tikzstyle{cblue}=[circle, draw, thin,fill=blue!20, scale=0.5]

\tikzstyle{new style 0}=[fill=black, draw=black, shape=circle]
\tikzstyle{red style 1}=[fill=red, draw=black, shape=circle]
\tikzstyle{blue style 2}=[fill=blue, draw=black, shape=circle]
\tikzstyle{white style 4}=[fill=white, draw=black, shape=circle]
\tikzstyle{bklack style 5}=[fill=black, draw=black, shape=rectangle]
\tikzstyle{red style 3}=[fill=red, draw=black, shape=rectangle]
\tikzstyle{yellow style 7}=[fill=yellow, draw=black, shape=rectangle]
\tikzstyle{new style 8}=[fill={rgb,255: red,0; green,132; blue,0}, draw={rgb,255: red,0; green,131; blue,0}, shape=circle]

\tikzstyle{new edge style 0}=[-]
\tikzstyle{new edge style 1}=[-, draw=red]
\tikzstyle{new edge style 2}=[-, draw=blue]
\tikzstyle{new edge style 3}=[-, draw={rgb,255: red,0; green,156; blue,0}]

\theoremstyle{definition}
\theoremstyle{plain}
\newtheorem{theorem}{\bf Theorem} 
\newtheorem{proposition}{\bf Proposition} 
\newtheorem{lemma}{\bf Lemma}
\newtheorem{corollary}{\bf Corollary} 
\newtheorem{conjecture}{\bf Conjecture}
\newtheorem{problem}{\bf Problem}

\newcommand\red[1] {{\color{red} #1}}
\newcommand\green[1] {{\color{green} #1}}
\newcommand\blue[1] {{\bf \color{blue} #1}}

\usepackage{stmaryrd}

\newcommand \brk[1]
{
	\llbracket #1\rrbracket
}

\def \R {{\mathbb R}}

\def \usecolour
{\newcommand {\rered}{\red}
\newcommand {\reblue} {\blue}
\newcommand {\regreen} {\green}
}

\def \nocolour
{\renewcommand {\rered}{}
\renewcommand {\reblue}{}
\renewcommand {\regreen}{}
}
\usecolour

\newcommand \them[2]
{\begin{theorem}
\label{#1} #2
\end{theorem}
}

\newcommand \lemm[2]
{\begin{lemma}
\label{#1} #2
\end{lemma}
}

\newcommand \prop[2]
{\begin{proposition}
		\label{#1} #2
	\end{proposition}
}

\newcommand \cor[2]
{\begin{corollary}
\label{#1} #2
\end{corollary}
}

\newcommand \equ[2]
{\begin{equation}\label{#1}
		#2
	\end{equation}
}

\newcommand \aln[1]
{\begin{align*}
		#1
	\end{align*}
}

\newcommand \Floor[1]
{
\left	\lfloor {#1}\right \rfloor 
}

\newcounter{countcase}
\def\incase{\addtocounter{countcase}{1}
{\noindent {\bf Case \thecountcase}: }}

\newcounter{countclaim}
\def\inclaim{\addtocounter{countclaim}{1}
{\vspace{0.2 cm}\noindent {\bf Claim \thecountclaim}: }}

\newcommand{\myproof}{{\noindent {\em Proof}.\quad}\setcounter{countclaim}{0}\setcounter{countcase}{0}}

\newcommand{\proofend}{{\hfill$\Box$}}

 \nocolour

\def \iff {if and only if }

\def \N {{\mathbb N}^+}

\def \Q {{\mathbb Q}}
\def \R {{\mathbb R}^+}

\def \scre{{\mathscr W}}
\def \scrt{{\mathscr T}}

\newcommand{\A}{\mathcal{A}}

\newcommand \Bra[1]
{\left ( #1\right ) }

\usepackage[hidelinks]{hyperref}

\begin{document}

\baselineskip 0.65 cm

\title{Trees with Prescribed Maximum Degree and Spectral Radius
}

\author[1\thanks{Corresponding author. Email: fengming.dong@nie.edu.sg and donggraph@163.com.}]
{Fengming Dong}

\author[2\thanks{Email: ruixuezhang7@163.com.}]
{Ruixue Zhang}

\affil[1]{\footnotesize National Institute of Education,
	Nanyang Technological University, Singapore}

\affil[2]{\footnotesize School of Mathematics and Statistics, Qingdao University, China}
	
\date{}

\maketitle

\baselineskip 0.6 cm

\begin{abstract}
It is well known that the spectral radius $\rho(T)$ of a tree $T$ 
with at least $3$ vertices 
satisfies
$\frac 14\rho(T)^2+1<\Delta(T)\le \rho(T)^2$, where $\Delta(T)$ is the 
maximum degree of $T$. Let $\mathbb{P}$ denote the set of spectral radii of all non-trivial trees. We ask
whether, for every $\alpha\in \mathbb{P}\setminus \{1\}$ and every integer $r$ 
satisfying $\frac 14\alpha^2+1<r\le \alpha^2$,
there exists a tree $T$ such that  $\Delta(T)=r$ and $\rho(T)=\alpha$.

For a positive integer $r$ and a positive real number $\alpha$, let
$\scre_r(\alpha)$ be the smallest subset of $[0,\infty)$ containing
$\alpha$ and closed under the following operation: if
$q_1,\ldots,q_s$ are positive 
elements of
$\scre_r(\alpha)$, 
not necessarily distinct, 
and
$q=\alpha-\sum_{i=1}^{s}q_i^{-1}$, then $q\in\scre_r(\alpha)$ whenever
either $1\le s<r$ and $q\ge0$, or $s=r$ and $q=0$.
We prove that 
$0\in {\mathscr W}_r(\alpha)$
if and only if there exists a tree $T$ with 
$\Delta(T)\le r$ and  $\rho(T)=\alpha$. 
Moreover, there exists a tree $T$ with $\Delta(T)=r$ and
$\rho(T)=\alpha$ if and only if there are positive elements
$q_1,\ldots,q_r\in\mathscr W_r(\alpha)$, not necessarily
distinct, satisfying $\sum_{i=1}^{r}q_i^{-1}=\alpha$.
Consequently, 
$
\mathbb{P}=\{\alpha\geq1:0\in {\mathscr W}
_{\lfloor\alpha^2\rfloor}(\alpha)\}.
$
As an application,
we  show that  for integers 
$k,r\ge 2$,  there exists a tree $T$  with $\Delta(T)=r$ and $\rho(T)=\sqrt k$ if and only if $\frac 14 k+1<r\le k$.
\end{abstract}

\vskip 2mm



\vskip 2mm


\section{Introduction}

\indent All graphs in this article
are finite, simple and undirected.
For a graph $G$,
let $V(G)$, $E(G)$, 
$\Delta(G)$ and $A(G)$
denote its vertex set, 
edge set, maximum degree
and the adjacency matrix, 
 respectively.
If $|V(G)|=n$, 
the adjacency matrix 
$A(G)$ has $n$ real
eigenvalues, counted with multiplicity (see~\cite{spectra:2010}).
The largest eigenvalue of $G$, 
denoted by $\rho(G)$, 
is called the {\it spectral radius} of $G$.

The adjacency spectral radius is a basic graph invariant whose relation
to graph structure has been studied extensively 
(see  \cite{Brouwer, spectra:2010,Stevan2015}). 
It also arises in sufficient conditions for packing edge-disjoint
spanning trees \cite{Fan-D}, 
the analysis of complex networks \cite{Meg2015}, and the study
of mixing times for single-site dynamics on spin systems 
\cite{Hayes2006}.

In this article, we study adjacency spectral radii of trees.
Lov\'asz and Pelik\'an \cite{Lovasz:1973} determined the maximum and minimum spectral
radii among trees of a fixed order. Guo and Tan \cite{Guo2001}
 obtained an upper
bound in terms of order and matching number, while Simi\'c and
To\v{s}i\'c~\cite{Simi2005}
identified a tree of maximum spectral radius among
trees with prescribed order and maximum degree. The following bounds
relate the spectral radius of a tree to its maximum degree; a proof
can be found in Stevanovi\'c
\cite{Stevan2003}.

\begin{theorem}[\cite{Stevan2003}]
\label{known0}
For any tree $T$ with $|V(T)|\ge 3$, 
$	\frac 14 \rho(T)^2+1
<\Delta(T)\le \rho(T)^2.$
\end{theorem} 

\def \mbbp {\mathbb{P}}

A graph $G$ is said to be 
non-trivial if 
$|V(G)|\ge 2$. 
Let \(\mathbb P\) denote the set of values \(\rho(T)\)
attained by non-trivial trees \(T\).
Note that $1$ and $\sqrt 2$ are the two smallest members of
$\Rho$, which are 
$\rho(K_2)$ and $\rho(K_{1,2})$, respectively, and 
$\rho(T)>\sqrt 2$ holds for every tree $T$ with $|V(T)|\ge 4$.
In this article, we will study the following inverse problem:

\begin{problem}\label{conj1}
	Is it true that 
	for any 
	$\alpha\in \mathbb{P}\setminus \{1\}$ and any integer $r$
	satisfying the condition
	$\frac 14 \alpha^2+1<r\le \alpha^2$, 
	there exists a tree $T$ with 
	$\rho(T)=\alpha$ and 
	$\Delta(T)=r$?
\end{problem}

Problem~\ref{conj1} asks whether the necessary degree bounds in Theorem 1 are
sufficient for trees with a prescribed adjacency spectral radius.
This is an inverse existence problem for unweighted adjacency matrices.
It is related in spirit to the inverse eigenvalue problem of a graph,
which seeks to characterize the spectra of real symmetric matrices
whose nonzero off-diagonal entries occur exactly at the edges of a
fixed graph  (see \cite{Hogben2005}).

Let $\N$ and 
$\R$ denote the set of positive integers and 
the set of positive real numbers,
respectively,
 and write $\brk{s} =\{1,\ldots,s\}$ for
$s\in\N$.
For $r\in \N$ and $\alpha\in \R$,
put  $\scre^{(0)}_r(\alpha)
=\{\alpha\}$, and 
for each $n\in \N$, 
let $\scre^{(n)}_r(\alpha)
= \scre^{(n-1)}_r(\alpha)\cup 
\mathcal A_r^{(n)}(\alpha)$,
where $\mathcal A_r^{(n)}(\alpha)$ is the set of all
non-negative numbers of the form
$q=\alpha-\sum_{i=1}^{s}q_i^{-1}$,
where 
$s\in \brk{r}$,
$q_1,\ldots,q_s\in
\scre_r^{(n-1)}(\alpha)\setminus\{0\}$
are not necessarily distinct, and either $s<r$ or $q=0$.
Now let 
\equ{new-set2}
{
\scre_r(\alpha)=
\bigcup_{n\ge 0}
\scre^{(n)}_r(\alpha).
}


\prop{prop1-1}
{
For $r\in \N$ and $\alpha\in \R$,
$\scre_r(\alpha)$ is the smallest subset of $[0,\infty)$ containing
$\alpha$ and closed under the following operation: if
$q_1,\ldots,q_s$ are positive 
elements of 
$\scre_r(\alpha)$, 
not necessarily distinct, 
and
$q=\alpha-\sum_{i=1}^{s}q_i^{-1}$, then $q\in\scre_r(\alpha)$ whenever
either $1\le s<r$ and $q\ge0$, or $s=r$ and $q=0$.
}

\proof 
	For positive elements $q_1,\ldots,q_s\in
	\scre_r(\alpha)$,
	all these elements belong to
	$\scre_r^{(N)}(\alpha)$ for some $N\ge0$.
	Every admissible output therefore belongs to
	$\scre_r^{(N+1)}(\alpha)$, proving closure.
	
	If $S\subseteq[0,\infty)$ contains $\alpha$ and is closed
	under the same operation, induction on $n$ gives
	$\scre_r^{(n)}(\alpha)\subseteq S$ for every $n\ge0$.
	Taking the union yields $\scre_r(\alpha)\subseteq S$,
	proving minimality.
\proofend 

The following properties are immediate from the definition.

\prop{prop1-0}
{
	$\scre_1(1)=\{0,1\}$,
	$\scre_1(\alpha)=\{\alpha\}$
	for $\alpha\in \R\setminus \{1\}$, and 
	$\scre_r(\alpha)=\{\alpha\}$
	for $r\in \N$ and 
	$0<\alpha<1$.
}

Our first goal is to characterize, using $\scre_r(\alpha)$, the existence
of a tree with spectral radius $\alpha$ and maximum degree exactly $r$,
and the existence of such a tree with maximum degree at most $r$.

\them{main-th01}
{
	For any $\alpha\in \R$  and 
	$r\in \N$, 
	there exists a tree $T$ 
	with $\Delta(T)=r$ and $\rho(T)=\alpha$
	if and only if 
	there exist positive 
	elements  $q_1,q_2,\cdots,q_r$ 
	of $\scre_r(\alpha)$,
	not necessarily distinct, 
	such that 
	$\sum\limits_{i=1}^r q_i^{-1}=\alpha$.
}

\them{main-th0}
{
For any 
$\alpha\in \R$ and $r\in \N$, 
there exists a tree $T$ with  $\Delta(T)\le r$ and $\rho(T)=\alpha$
\iff $0\in \scre_r(\alpha)$.
}

\cor{main-th0-1}
{$\Rho=
\big \{ \alpha\ge 1: 
	0\in \scre_{\Floor{\alpha^2}}
	(\alpha)
	\big \}$.
}

\proof
If $\alpha=1$, then $\alpha=\rho(K_2)$ and $0\in \scre_1(1)$. 
If $\alpha>1$ belongs
to $\Rho$, any realizing tree has at least three vertices, and Theorem~\ref{known0}
gives
$\Delta(T)\leq\lfloor\alpha^2\rfloor$. 
Theorem~\ref{main-th0} and the inclusion
$\scre_{\Delta(T)}(\alpha)\subseteq \scre_{\lfloor\alpha^2\rfloor}(\alpha)$ give the
forward implication. Conversely, 
Theorem~\ref{main-th0} gives a realizing tree whenever
$\alpha\geq1$ and 
$0\in \scre_{\lfloor\alpha^2\rfloor}(\alpha)$.
\proofend 

For $\alpha>0$, let $\scre(\alpha)$ be the smallest subset of
$[0,\infty)$ containing $\alpha$ with the following property: for every
$s\in\N$ and all positive $q_1,\ldots,q_s\in\scre(\alpha)$, not
necessarily distinct, if
\[
q:=\alpha-\sum_{i=1}^{s}q_i^{-1}\geq0,
\]
then $q\in\scre(\alpha)$.
Equivalently, 
	$\scre(\alpha)
	=\bigcup\limits_{r\ge1}\scre_r(\alpha)$.
	
	Every element of $\scre(\alpha)$ is obtained by finitely many applications
	of the recursive rule. Hence it belongs to $\scre_r(\alpha)$ for some
	positive integer $r$.

\cor{main-th0-10}
{$\Rho=\big \{ \alpha\ge 1: 0\in \scre(\alpha)\big \}$.
}

\myproof 
 If $\alpha\in \Rho$, 
then $\alpha\ge 1$ and 
$0\in \scre_{\Floor{\alpha^2}}
(\alpha)
\subseteq \scre(\alpha)$.
Conversely, if 
$0\in\scre(\alpha)$, 
then  there exists $r\in \N$ 
such that 
$0\in \scre_r(\alpha)$,
and Theorem~\ref{main-th0}
gives $\alpha\in \Rho$.
\proofend 

Our second aim is to use the recursive characterization to answer
Problem~\ref{conj1} when $\alpha=\sqrt{k}$ for an integer $k\geq2$.
In particular, this covers all integer spectral radii $\alpha\geq2$.

\them{main-th2}
{
	For  integers $k,r\ge 2$, 
	there exists a tree $T$ with $\rho(T)=\sqrt k$ and $\Delta(T)=r$
	\iff $\frac 14 k+1<r\le k$.
}

\cor{main-th2-r1}
{
For integers $k, r\ge 2$, 
there exists a tree $T$ with $\rho(T)=k$ and $\Delta(T)=r$
\iff $\frac 14 k^2+1<r\le k^2$.
}

The remainder of this article is structured as follows. 
Section~\ref{sec2}
proves the necessity of 
Theorem~\ref{main-th01}.
Section~\ref{sec3}
proves the sufficiency in 
Theorem~\ref{main-th01}
and establishes
Theorem~\ref{main-th0}, and 
Section~\ref{sec5} proves Theorem~\ref{main-th2}.
Finally, 
Section~\ref{problem} 
establishes further results on sums and averages of square roots
and poses two related open problems.

\section{The necessity of Theorem~\ref{main-th01}
\label{sec2}
} 

\subsection{A condition 
	guaranteeing $\rho(G)=\alpha$}

The Perron–Frobenius theorem, proved by Oskar Perron \cite{Perron:1907}
in 1907 and Georg Frobenius \cite{Frob:1912} in 1912,
is an important result 
in the study of eigenvalues 
and eigenvectors of square matrices. 
Applying the 
Perron–Frobenius theorem for irreducible non-negative matrices, 
the following conclusion for the spectral radius $\rho(G)$ of a connected 
graph $G$ is obtained.
Let $\lambda_i(G)$ denote the $i$-th 
largest eigenvalue of $G$. 
Then, $\rho(G)=\lambda_1(G)$.

\begin{theorem}[Perron–Frobenius theorem] 
	\label{PT}
	Let $G$ be a connected non-trivial graph. Then  $\lambda_1(G)>\max\{0,\lambda_2(G)\}$ and  there is an eigenvector $\bar x$ of $A(G)$ corresponding to $\lambda_1(G)$ in which all 
	components are positive. 
	Furthermore, $\lambda_1(G)$ is the only eigenvalue of $G$ admitting an eigenvector 
	all of whose components are positive.
\end{theorem}

Applying Theorem~\ref{PT}, 
the next conclusion follows
directly.

\prop{prop6-0}
{Let $G$ be a connected non-trivial 
	graph
	and $\alpha$ be a positive number.
	Then,  $\rho(G)=\alpha$
	\iff  there is a mapping 
	$\phi:V(G)\rightarrow \R$  such that  
	for each $u\in V(G)$, 
	\equ{p6-0-e1}
	{	
	\alpha \phi(u)
	=\sum_{w\in N(u)} \phi(w).
}
}

\subsection{
Necessity in Theorem~\ref{main-th01}
}

For a fixed $\alpha$, equation~(\ref{p6-0-e1}) is the eigenvector equation
$A(G)\phi=\alpha\phi$.
If $G$ is connected and $\phi$ is positive, then
Proposition~\ref{prop6-0} gives $\alpha=\rho(G)$.

For a directed tree $\vect$
and $v\in V(\vect)$,
let $id_{\vect}(v)$ denote 
the in-degree 
of $v$ in $\vect$
(i.e., the number of arcs 
whose head is $v$).
If $id_{\vect}(v)=0$, then $v$ is called 
a source of $\vect$.
Let $\sett$ be the set of directed trees $\vect$ 
that have exactly one source.
Clearly, for any $\vect\in \sett$, 
if $u$ is the unique source of $\vect$, 
then $id_{\vect}(v)=1$ 
for each $v\in V(\vect)\setminus \{u\}$.
Note that 
for any undirected tree 
$T$ and any vertex $u$ in this tree, orienting every edge away from \(u\) is the unique
orientation with a sole 
source \(u\).
We now 
prove the following result
from which 
the necessity of Theorem~\ref{main-th01}
follows directly.

\prop{pro2-1}
{
	Let $T$ be a  tree  with $|V(T)|\ge 2$, 
	$\alpha=\rho(T)$
	and $r=\Delta(T)$, and let
	$\phi: V(T)\rightarrow \R$ 
	be a positive eigenvector of $A(T)$ 
	with respect to $\alpha$. 
	Then the following conclusions hold:
	\begin{enumerate}
				[itemsep=0mm, label=(\roman*)]
		\item 
		$\frac{\phi(v')}{\phi(v)}\in \scre_r(\alpha)$
		for every ordered pair 
		$(v',v)$ of 
		adjacent vertices
		 in $T$; 
		\item $\alpha^{-1}\in \scre_r(\alpha)$; 
		and 
		\item if $s=d_T(u)$ for some vertex $u$ in $T$,  
		then there exist
		positive elements 
		$q_1,\ldots,q_s$ of 
		 $\scre_r(\alpha)$,
		 not necessarily distinct, 
		such that $\sum\limits_{i=1}^s q_i^{-1}
		=\alpha$.
	\end{enumerate}
}

\myproof We will first 
establish two claims.
Let $\vect$ be any orientation of $T$ with a unique source.
Thus $\vect\in \sett$.
Clearly, for any $v\in V(\vect)$, 
$id_{\vect}(v)=0$
if $v$ is the source of $\vect$,
and 
$id_{\vect}(v)=1$ otherwise.
For each $v\in V(\vect)$ with $id_{\vect}(v)=1$, let 
$\Psi(v)\coloneq\frac{\phi(v')}{\phi(v)}$,
where $v'$ is the only vertex in $\vect$ such that $(v',v)$ 
is a directed edge in $\vect$.

For each vertex in $\vect$,
let $od_{\vect}(v)$ denote the 
out-degree of $v$ in $\vect$.

\inclaim $\Psi(v)=\alpha$ 
for each $v\in V(\vect)$
with $od_{\vect}(v)=0$.

Since $od_{\vect}(v)=0$, 
$v'$ is the unique neighbor of $v$ in $T$. 
By (\ref{p6-0-e1}),
\equ{pro2.6-e2}
{
	\alpha\, \phi(v)=\phi(v'),
}
implying that $\Psi(v)=\frac{\phi(v')}{\phi(v)}
=\alpha$.
Claim 1 holds.

\inclaim $\Psi(v)\in \scre_r(\alpha)$ for each 
$v\in V(T)$ with $id_{\vect}(v)=1$.

We prove Claim 2 by induction on the height of a non-source vertex, where its height is the maximum length of a directed path starting at that vertex.
First, 
by Claim 1, $\Psi(v)=\alpha\in \scre_r(\alpha)$  if $od_{\vect}(v)=0$.
Now assume that $v$ is a vertex 
in $T$  with 
$id_{\vect}(v)=1$ and 
$N^+_{\vect}(v)
=\{v_1,\cdots,v_t\}$,
where $t\coloneq od_{\vect}(v)>0$
and $N^+_{\vect}(v)$ is the 
set of vertices $w$ in $\vect$ 
with $(v,w)\in E(\vect)$.
Assume that 
$\Psi(v_i)\in \scre_r(\alpha)$ for each $i\in \brk{t}$.
By the definition of $\Psi$, 
$\Psi(v_i)=\frac{\phi(v)}{\phi(v_i)}$.
By (\ref{p6-0-e1}), 
\begin{equation}
	\label{pro2.6-e3}
	\alpha\, \phi(v)
	=\phi(v')+
	\sum_{i=1}^t \phi(v_i).
\end{equation}
It follows that 
\equ{pro2.6-e4}
{
0<	\Psi(v)=\frac{\phi(v')}{\phi(v)}
	=\alpha-\sum_{i=1}^t
	\frac{\phi(v_i)}{\phi(v)}
	=\alpha-\sum_{i=1}^t
	\Psi(v_i)^{-1}.
}
As $t=od_{\vect}(v)=d_T(v)-1<r$, 
by the definition of $\scre_r(\alpha)$, $\Psi(v)\in \scre_r(\alpha)$.
Hence Claim 2 holds.

We now prove (i), (ii) and (iii).

For any ordered pair $(v',v)$ 
of adjacent vertices $v'$ and $v$ in $T$, there exists 
an orientation $\vect$ with exactly one source such that 
$(v',v)$ is a directed edge in $\vect$.
By Claim 2, $\frac{\phi(v')}{\phi(v)}=\Psi(v)\in \scre_r(\alpha)$.
Thus, (i) holds.

Choose a leaf $v_1$ of $T$, and let $v_2$ be its unique
neighbour. 
By equation~(\ref{p6-0-e1}),
$\alpha\phi(v_1)=\phi(v_2)$, so
$\phi(v_1)/\phi(v_2)
=\alpha^{-1}$.
Part~(i) therefore gives
$\alpha^{-1}\in\scre_r(\alpha)$.
Hence (ii) holds.

Assume that $u$ is a vertex in $T$ and $s=d_T(u)$.
Let 
$\vect$ be the orientation of $T$ 
such that $u$ is its unique source. Let
$N^+_{\vect}(u)
=\{u_1,\cdots,u_s\}$.
By Claim~2, for each $i\in\brk{s}$, the number
$
q_i:=\frac{\phi(u)}{\phi(u_i)}=\Psi(u_i)
$
is a positive element of $W_r(\alpha)$.
By (\ref{p6-0-e1}), 
\begin{equation}
	\label{pro2.6-e5}
	\alpha\, \phi(u)
	=\sum_{i=1}^s \phi(u_i).
\end{equation}
It follows that 
\equ{}
{
	\sum_{i=1}^{s}q_i^{-1}
	=\sum_{i=1}^s \frac{\phi(u_i)}{\phi(u)}
	=\alpha.
}
Thus conclusion (iii) holds.
\proofend 

\vspace{3 mm}

By taking $u$ to be a vertex in $T$ with $d_T(u)=\Delta(T)$, 
Proposition~\ref{pro2-1}
yields 
the following conclusion,
and thus
the necessity part of
Theorem~\ref{main-th01} 
follows directly.

\cor{se2-co1}
{
	Let $T$ be a tree of order at least $2$, 
	and let $\alpha=\rho(T)$
	and $r=\Delta(T)$.
	Then 
	there exist
	positive elements 
	$q_1,\ldots,q_r$ of 
	$\scre_r(\alpha)$,
	not necessarily distinct, 
	such that $\sum\limits_{i=1}^r q_i^{-1}
	=\alpha$.
}

\section{Proofs of Theorems~\ref{main-th01}
	and~\ref{main-th0}
\label{sec3}
}

In this section, we always assume that  $\alpha\in \R$ 
and $r\in \N$.

\subsection{Directed trees $\vect$ with a mapping $\omega: V(\vect)\to \scre_r(\alpha)$
}

Let $\scrt_{r,\alpha}$ denote the
set of 
ordered pairs $(\vect,\omega)$,
where $\vect\in \sett$ is a directed tree 
with a unique source $z$
and $\omega$ is a 
mapping from $V(\vect)$ to 
$\scre_r(\alpha)$
such that for each $v\in V(\vect)$, 
\begin{enumerate}
		[itemsep=0mm, label=(\roman*)]
\item 
$od_{\vect}(v)\le r$ if $v=z$ and $\omega(z)=0$, 
and $od_{\vect}(v)\le r-1$ otherwise; 
and 
\item 
$\omega(v)>0$ whenever 
$v\in V(\vect)\setminus \{z\}$, and 
\equ{map0}
{
	\omega(v)=
	\alpha-
	\sum\limits_{u\in N^+_{\vect}(v)}
	\omega(u)^{-1}.
}
\end{enumerate} 
Obviously, if $od_{\vect}(v)=0$, then $N^+_{\vect}(v)=\emptyset$, 
implying that $\omega(v)=\alpha$.

For $\beta\ge 0$, let 
$\scrt_{r,\alpha}(\beta)$ denote 
the set of ordered pairs $(\vect,\omega)$ in $\scrt_{r,\alpha}$ 
such that $\omega(z)=\beta$,
where $z$ is the unique source of $\vect$.
By the definitions of  $\scre_{r}(\alpha)$ and 
$\scrt_{r,\alpha}(\beta)$, 
the following conclusion holds.

\prop{lem-n1}
{
For any $\beta\in \scre_{r}(\alpha)$, 
$\scrt_{r,\alpha}(\beta)\ne \emptyset$.
}
 
\myproof By definition, 
$\scre_{r}(\alpha)=
\bigcup_{n\ge 0}
\scre^{(n)}_{r}(\alpha)$.
Thus, for each $\beta\in
\scre_{r}(\alpha)$, we have 
$\beta\in 
\scre^{(n)}_{r}(\alpha)$
for some $n\ge 0$. 
The proof will be completed
by induction on $n$.
If $\beta\in  \scre^{(0)}_{r}(\alpha)$, 
then $\beta=\alpha$
and  $(\vect_0,\omega)\in  \scrt_{r,\alpha}(\beta)$,
where $\vect_0$ is the one-vertex tree and $\omega(z)=\alpha$ for 
the unique vertex $z$ in $\vect_0$.

Assume that 
the result holds for 
each  $\beta\in \scre^{(n)}_{r}(\alpha)$, 
where $n\ge 0$.
If $\scre^{(n+1)}_{r}(\alpha)
=\scre^{(n)}_{r}(\alpha)$,
 the induction
step is immediate.
Assume that 
$\scre^{(n+1)}_{r}(\alpha)
\ne \scre^{(n)}_{r}(\alpha)$
and $\beta\in 
\scre^{(n+1)}_{r}(\alpha)
\setminus 
\scre^{(n)}_{r}(\alpha)$.
By definition, 
there exist $s\in \brk{r}$ and 
positive elements $\beta_1,
\ldots, \beta_s$ of 
$\scre_r^{(n)}
(\alpha)$,
not necessarily distinct, 
such that
\[
\beta=\alpha-\sum_{i=1}^{s}\beta_i^{-1}\geq0,
\]
where $\beta=0$ whenever $s=r$.
 By the inductive assumption, 
 for each $i\in \brk{s}$, 
 there exists $(\vect_i, \omega_i) \in \scrt_{r,\alpha}(\beta_i)$.
 Take pairwise
 vertex-disjoint copies of these trees 
 $\vect_1, \cdots, \vect_s$
 and 
let $\vect$ be the directed tree
obtained 
from them
by adding a new vertex $z$ and a new directed edge 
$(z,z_i)$ for each $i\in \brk{s}$, where $z_i$ is the unique source of $\vect_i$,
and let $\omega$ be 
the mapping from $V(\vect)$ to 
$\scre_{r}(\alpha)$ such that 
$\omega(v)=\omega_i(v)$ holds
for each $i\in \brk{s}$ and 
$v\in V(\vect_i)$, 
and
\equ{define-w}
{
\omega(z)
=
\alpha-\sum\limits_{i=1}^s \omega(z_i)^{-1}
=\alpha-\sum\limits_{i=1}^s \beta_i^{-1}=\beta.
}
Since each old root $z_i$ 
has a positive label $\beta_i$, 
it has at most $r-1$ children; after the new edge to its
parent is added, 
its degree is at most $r$. 
All other old vertices retain degree at most $r$, and the
new root has degree $s \le r$.
Obviously, $(\vect,\omega)\in \scrt_{r,\alpha}(\beta)$.
Hence the proposition holds.
\proofend 

For example, two members 
$(\vect_i,\omega_i)$ for $i=1,2$
of  $\scrt_{3,2}(\frac{1}{2})$
are shown in Figure~\ref{F3-1},
where the value of $\omega_i(v)$ 
for each vertex $v$ in $\vect_i$ is shown beside vertex $v$.

\begin{figure}[H]
	\centering
	\begin{tikzpicture}
		\tikzset{vertex/.style = {shape=circle,fill, draw,minimum size=1em}}
		\tikzset{edge/.style = {->,> = latex'}}

		\tikzstyle{vertex}=[circle, fill=black, inner sep=0pt, minimum size=8pt]

		\node[vertex, label=left:$\frac{1}{2}$] (a) at (0,0) {};
		\node[vertex, label=above:1] (b) at (1.6,0.6) {};
		\node[vertex, label=right:2] (c) at (1.6,-0.6) {};
		\node[vertex, label=right:2] (d) at (3.2,1.2) {};
		\node[vertex, label=right:2] (e) at (3.2,0) {};
		
		\draw[edge, thick] (a) to (b);
		\draw[edge, thick] (b) to (d);
		\draw[edge, thick] (b) -- (e);
		\draw[edge, thick] (a) -- (c);
		
	\end{tikzpicture}
	\hspace{0.1\textwidth}
	\centering
	\begin{tikzpicture}
		\tikzset{vertex/.style = {shape=circle,fill, draw,minimum size=1em}}
		\tikzset{edge/.style = {->,> = latex'}}
		
		\tikzstyle{vertex}=[circle, fill=black, inner sep=0pt, minimum size=8pt]
		
		\node[vertex, label=above:$\frac{1}{2}$] (A) 
		at (0,0) {};
		\node[vertex, label=above:$\frac{2}{3}$] (B) at (1.6,0) {};
		\node[vertex, label=above:$\frac{3}{2}$] (C) at (3.2,0.7) {};
		\node[vertex, label=right:$2$] (D) at (4.8,0.7) {};
		\node[vertex, label=below:$\frac{3}{2}$] (E) at (3.2,-0.7) {};
		\node[vertex, label=right  :$2$] (F) at (4.8,-0.7) {};
		
		\draw[edge, thick] (A) -- (B);
		\draw[edge, thick] (B) -- (C);
		\draw[edge, thick] (C) -- (D);
		\draw[edge, thick] (B) -- (E);
		\draw[edge, thick] (E) -- (F);
		
	\end{tikzpicture}
	
	\centerline{~~~~(a) $(\vect_1,\omega_1)$ \hspace{4 cm} (b) $(\vect_2,\omega_2)$  ~~~~ }
	
	\caption{$(\vect_i,\omega_i)
		\in \scrt_{3,2}(\frac{1}{2})$
		for $i=1,2$}
	\label{F3-1}
\end{figure}

By Proposition~\ref{lem-n1}, if $0\in \scre_r(\alpha)$, then 
$\scrt_{r,\alpha}(0)\ne \emptyset$,
and there exists $(\vect,\omega)$ 
in $\scrt_{r,\alpha}(0)$.
For any directed tree $\vect$,
write $\rho(\vect)=\rho(T)$
and $\Delta(\vect)=\Delta(T)$,
where $T$ is the underlying 
graph of $\vect$.

\prop{lem-n2}
{
	For any $(\vect,\omega)\in \scrt_{r,\alpha}(0)$,
	we have 
	$\Delta(\vect)\le r$ and 
	$\rho(\vect)=\alpha$.
}

\myproof By the definition
of $\scrt_{r,\alpha}(0)$, 
$\Delta(\vect)\le r$. 
It suffices to show that $\rho(\vect)=\alpha$. 

By definition again, $\vect$ is a directed 
tree with a unique source $z$ 
and 
$\omega(v)\in \scre_{r}(\alpha)$
for each $v\in V(\vect)$.
If \(od_{\vect}(z)=0\), 
equation~\eqref{map0} gives
\(\omega(z)=\alpha>0\), contrary to \(\omega(z)=0\).
Thus, $s\coloneq od_{\vect}(z)\ge 1$.
Let $\vect_1,\vect_2,\dots, \vect_s$ be the components of $\vect-z$.
Clearly, $1\le s\le r$.
For each $i\in \brk{s}$, 
$(\vect_i, \omega_i)\in 
\scrt_{r,\alpha}(\beta_i)$, where 
$\omega_i$ is the restriction 
of $\omega$ to  $V(\vect_i)$, 
$\beta_i=\omega(z_i)$
and $z_i$ is the vertex in $\vect_i$ 
with $z_i\in N^+_{\vect}(z)$.
Then $z_i$ is the unique source of $\vect_i$ for each $i\in \brk{s}$ and 
\equ{lem-n2-e1}
{
	0=\omega(z)
	=\alpha-\sum_{i=1}^s \omega(z_i)^{-1}
	=\alpha-\sum_{i=1}^s \beta_i^{-1}.
}

Now let $\phi$ be the mapping from 
$V(\vect)$ to $\R$ defined below: 
\vspace{-2 mm} 
\begin{enumerate}
			[itemsep=0mm, label=(\roman*)]
	\item $\phi(z)=1$; and 	
	\item  
	for each directed edge 
	$(w_1,w_2)\in E(\vect)$,
	$\phi(w_2)=\frac{\phi(w_1)}
	{\omega(w_2)}$.
\end{enumerate}
Every non-source vertex has a unique parent, so  $\phi$
is uniquely defined 
by induction on the distance from the source. All its values are positive because every non-source label $\omega(v)$ is positive.

In order to show that $\rho(\vect)=\alpha$,
by Proposition~\ref{prop6-0}, 
it suffices to prove the following claim. 

\noindent {\bf Claim A}: 
$\phi$ 
satisfies equation  (\ref{p6-0-e1}) for each vertex $u\in V(\vect)$.

We first prove Claim A 
for the case $u=z$ or $u$ is a 
sink of $\vect$ (i.e., $N^+_{\vect}(u)=\emptyset$).
As $N^+_{\vect}(z)=\{z_i: i\in \brk{s}\}$, 
we have 
\equ{lem-n2-e10}
{
	\sum_{v\in N^+_{\vect}(z)}\phi(v)=
	\sum_{i=1}^s\phi(z_i)=
	\sum_{i=1}^s\frac{\phi(z)}{\omega(z_i)}
	= \phi(z)\sum_{i=1}^s\beta_i^{-1}
	=\alpha
	\phi(z) ,
}
where the last step follows from 
(\ref{lem-n2-e1}).

For the case that $u$ is a sink 
of $\vect$,  $\omega(u)=\alpha$ by
(\ref{map0}).
If $u'$ is the unique in-neighbor 
of $u$ in $\vect$, 
then 
$\phi(u')=\phi(u)\omega(u)$ by the definition of $\phi$, and 
\equ{lem-n2-e2}
{
	\sum_{v\in N_{T}(u)}\phi(v)=
	\phi(u')=
	\phi(u)\omega(u)
	=\alpha \phi(u).
}
Finally, we assume that $u$ is a vertex in $\vect$ which is neither the source 
nor a sink of $\vect$.
Assume that $u'$ is the only in-neighbor of $u$ in $\vect$ 
and 
$N^+_{\vect}(u)=\{u_i:i\in \brk{t}\}$.
By the definition of $\phi$, 
\equ{p7-0-e3}
{
	\phi(u')+\sum_{i=1}^t\phi(u_i)
	=\omega(u)\phi(u)
	+\sum_{i=1}^t \frac{\phi(u)}{\omega(u_i)}
	=\phi(u) 
	\left ( 
	\omega(u)+\sum_{i=1}^t \omega(u_i)^{-1}
	\right )
	=\alpha\phi(u),
}
where the last step follows from
(\ref{map0}).
Thus, $\phi$ satisfies 
equation (\ref{p6-0-e1})  
for each $u\in V(\vect)$, 
implying that 
$\rho(\vect)=\alpha$ by Proposition~\ref{prop6-0}.
\proofend 

\subsection{Proofs of 
	Theorems~\ref{main-th01}
	and~\ref{main-th0}}

In this subsection, we will complete the proofs of 
Theorems~\ref{main-th01}
and~\ref{main-th0}.
We first find a sufficient condition for the existence 
of a tree  $T$ 
with $\Delta(T)\le r$ and $\rho(T)=\alpha$.

\prop{corn-1}
{If there are positive 
	elements $q_1,q_2,\ldots,q_s$
	of $\scre_r(\alpha)$, 
	not necessarily distinct, 
	where $s\in \brk{r}$,
	such that 
	$\sum\limits_{i=1}^sq_i^{-1}=\alpha$, 		
	then there exists a tree $T$ 
	with $s\le \Delta(T)\le r$ and $\rho(T)=\alpha$.
}

\myproof Assume that 
$\sum\limits_{1\le i\le s}q_i^{-1}=\alpha$, where 
$q_i\in \scre_r(\alpha)\setminus\{0\}$ 
for each $ i\in \brk{s}$.
By Proposition~\ref{lem-n1}, 
for each $i\in \brk{s}$,
there exists $(\vect_i,\omega_i)\in \scrt_{r,\alpha}(q_i)$.
Assume that $z_i$ is the only source of $\vect_i$ for each $i\in \brk{s}$, and let 
$\vect$ be the directed tree 
obtained from $\vect_1,\cdots,\vect_s$,
which are pairwise vertex-disjoint, 
by adding a new vertex $z$ 
and directed edges $(z,z_i)$ 
for each $i\in \brk{s}$, and 
let $\omega$ be the mapping 
from $V(\vect)$ to $\scre_r(\alpha)$ 
defined as follows: 
$\omega(z)=0$ and 
$\omega(v)=\omega_i(v)$ 
whenever $v\in V(\vect_i)$,
where $i\in \brk{s}$.
By the definition of $\vect$ and $\omega$, 
$(\vect,\omega)\in \scrt_{r,\alpha}(0)$.
Clearly, 
$r\ge \Delta(\vect)\ge od_{\vect}(z)=s$.
By Proposition~\ref{lem-n2}, $\rho(\vect)=\alpha$.
Thus, the result holds
for the underlying tree $T$ of $\vect$. 
\proofend 

The following corollary shows 
some statements equivalent to 
the existence of a tree $T$ with $\Delta(T)\leq r$ and $\rho(T)=\alpha$,
from which Theorem~\ref{main-th0} follows.

\begin{corollary}
	\label{cor3-1}
	For $r\in \N$ and $\alpha\ge 1$, 
	the following statements are equivalent:
	\begin{enumerate}
				[itemsep=-1mm, label=(\roman*)]
		\item there exists a tree $T$ with $\Delta(T)\leq r$ and $\rho(T)=\alpha$;
		\item $\alpha^{-1}\in \scre_r(\alpha)$;
				\item	
				 there exist $a,b\in \scre_r(\alpha)$ with $ab=1$;
		\item $0\in \scre_r(\alpha)$;
		and 
		
		\item 
			there are positive 
		elements $q_1,q_2,\ldots,q_s$
		of $\scre_r(\alpha)$, 
		not necessarily distinct, 
		where $s\in \brk{r}$, such that $\sum\limits_{i=1}^{s}q_i^{-1}=\alpha$.
	\end{enumerate}
\end{corollary}

\proof First suppose that 
$r = 1$. 
By Proposition~\ref{prop1-0}, 
$\scre_1(1)=\{0,1\}$ 
and $\scre_1(\alpha)=\{\alpha\}$
for any $\alpha\in \R\setminus \{1\}$.
Thus,
each of the statements (i) - (v) holds, and 
the conclusion follows.
In the remainder of the proof, 
assume that $r\ge 2$.

(i) $\Rightarrow$ (ii).
Let $T$ be a tree as in (i), and put $r':=\Delta(T)\le r$.
By  Proposition~\ref{pro2-1},
$\alpha^{-1}\in 
\scre_{r'}(\alpha)
\subseteq \scre_{r}(\alpha)$.

(ii) $\Rightarrow$ (iii).  Since $\alpha^{-1}\in \scre_r(\alpha)$,
(iii) holds by taking $a=\alpha$ and $b=\alpha^{-1}$.

(iii) $\Rightarrow$ (iv). 
If $\alpha=1$, then $0=1-\frac 11\in \scre_r(\alpha)$.

For $\alpha>1$, exchange $a,b$ if necessary so that $0<b\le1<\alpha$.
Since 
$b\ne \alpha$, 
let $n$ be the least integer 
for which $b\in \scre^{(n)}_r(\alpha)$.
By the recursive
definition, there exist positive elements
$q_1,\ldots,q_t\in  \scre^{(n-1)}_r(\alpha)\subseteq \scre_r(\alpha)$
such that 
$$
b=\alpha-\sum_{i=1}^t \frac 1{q_i}.
$$
As $b>0$, 
the admissibility condition gives $1\le t<r$.
Now $b=a^{-1}$ gives 
$$
0=\alpha-a^{-1}-\sum_{i=1}^{t}q_i^{-1}.
$$
There are $t + 1 \le r$ positive inputs, so this zero output is admissible. 
Since $a\in \scre_r(\alpha)$
and $q_i\in \scre_r(\alpha)$
for each $i\in \brk{t}$, 
by Proposition~\ref{prop1-1},
the above expression implies 
that 
$0\in \scre_r(\alpha)$.
This proves (iv).

(iv) $\Rightarrow$ (v). 
Since $\alpha>0$, zero is not in $\mathscr W_r^{(0)}(\alpha)$. Let $n\ge1$ be the least integer such that $0\in\mathscr W_r^{(n)}(\alpha)$. The recursive definition then gives $1\le s\le r$ and positive elements $q_1,\ldots,q_s\in\mathscr W_r^{(n-1)}(\alpha)$ satisfying
\[
0=\alpha-\sum_{i=1}^{s}q_i^{-1}.
\]
Thus statement~\textup{(v)} holds.

(v) $\Rightarrow$ (i). 
It follows from Proposition~\ref{corn-1}.
\proofend

We can now complete 
the proofs of Theorems~\ref{main-th01}
and~\ref{main-th0}.

\vspace{3 mm}

\noindent {\it Proofs of Theorems~\ref{main-th01}
and~\ref{main-th0}}:
Corollary~\ref{se2-co1} proves the necessity in
Theorem~\ref{main-th01}.
By Proposition~\ref{corn-1},
its sufficiency also holds
by taking $s$ to be $r$.

If $0<\alpha<1$, then Proposition~\ref{prop1-0} gives $\scre_r(\alpha)=\{\alpha\}$, so
$0\notin \scre_r(\alpha)$. Every tree of order at least two has spectral radius at
least $1$, whereas the one-vertex tree has spectral radius $0$. Thus both sides
of Theorem 3 are false for $0<\alpha<1$. For $\alpha\geq1$, the theorem follows
from Corollary~\ref{cor3-1}.
\proofend

\section{Proof of Theorem~\ref{main-th2}
\label{sec5}
} 

Let $k\in \N$ with $k\ge 2$
and $r_0=\Floor{\frac 14 k}+2$.
Clearly, for any $r\in \N$,
$r>\frac 14 k+1$ \iff $r\ge r_0$.
Now assume that $k=4m+t$, where $m\in \{0\}\cup \N$ and 
$t\in \{0,1,2,3\}$.
Then $r_0=\Floor{\frac {1}{4}k}+2=m+2$.
Let 
\equ{se4-e1}
{
\A=\left\{x\ge 0:\frac{x}{\sqrt{k}}
\in\scre_{m+2}(\sqrt{k})\right\}
}
and for $x_1,\ldots,x_\ell\in
\A\setminus \{0\}$, 
not necessarily distinct, 
let
\equ{se4-e1}
{
\Gamma(x_1,\ldots,x_\ell)
=k\Bra{1-\sum_{j=1}^{\ell}\frac{1}{x_j}}.
}
We set $\Gamma() = k$ when no arguments are present.

\lemm{se4-l1}
{We have $k\in \A$.
	Moreover, 
	if $x_1,x_2,\cdots,x_{\ell}\in \A\setminus \{0\}$,
	$\ell\in \brk{m+1}$, and 
	$\Gamma(x_1,\ldots,x_\ell)
	\ge 0$, then
$
\Gamma(x_1,\ldots,x_\ell)\in\A.     $
}

\myproof 
As $\sqrt k\in \scre_{m+2}(\sqrt k)$, we have $k\in \A$.
Since $x_i\in \A\setminus \{0\}$,
$\frac{x_i}{\sqrt k}\in \scre_{m+2}(\sqrt k)$
for each $i\in \brk{\ell}$.
By the definition of 
$\scre_{m+2}(\sqrt k)$,
if $\ell\in \brk{m+1}$
and $\Gamma(x_1,\ldots,x_\ell)
\ge 0$,
then 
$$
0\le 
\frac{\Gamma(x_1,\ldots,x_\ell)}{\sqrt k}
=\sqrt k- \sum_{j=1}^{\ell} \Bra{x_j/\sqrt k}^{-1} 
\in \scre_{m+2}(\sqrt k).
$$
Thus, the lemma holds.
\endproof

For any $x\in \A$, let $x^{[q]}$ 
denote $q$ copies of $x$.
Thus, by Lemma~\ref{se4-l1}, 
 for $x_1,x_2,x_3\in \A\setminus \{0\}$
and  $q_1,q_2,q_3\in \{0\}\cup \N$ 
with $q_1+q_2+q_3\le m+1$, 
if $\Gamma(x_1^{[q_1]},x_2^{[q_2]}, x_3^{[q_3]})\ge 0$, then 
\equ{se4-e0}
{
\Gamma(x_1^{[q_1]},x_2^{[q_2]}, x_3^{[q_3]})
=k\Bra{1-\frac {q_1}{x_1}
-\frac {q_2}{x_2}
-\frac {q_3}{x_3}
}\in \A.
}

\prop{le8-1}
{The following conclusions hold:
\begin{enumerate} 
			[itemsep=0mm, label=(\roman*)]
	\item If $m=0$ and $2\le t\le 3$, 
	$s\in \A$ holds
	for every integer $s$ with
	$2\le s\le t$; and 
	\item if $m\ge 1$ and $t\in \{0,1,2,3\}$, 
	$s\in \A$ holds
	for every integer $s$ 
	with 
	$3m+t-1\le s\le 4m+t$.
\end{enumerate} 
}

\myproof (i). 
Assume that $m=0$.
By Lemma~\ref{se4-l1}, 
$k=t\in \A$. 
If $t=3$, then $k=3$ and by
(\ref{se4-e0}), 
$$
2=3\Bra{1-\frac 13}
=\Gamma(3)
\in \A.
$$

(ii) Assume that $m\ge 1$. 
By Lemma~\ref{se4-l1}, 
$k=4m+t\in \A$. 
For any $s\in \N$ with 
$3m+t-1\le s\le 4m+t-1$, 
we have 
$1\le 4m+t-s\le m+1$, 
and by (\ref{se4-e0}), 
$$
s=(4m+t)\Bra{1
	-\frac{4m+t-s}{4m+t}}
=
\Gamma(k^{[4m+t-s]})\in \A.
$$
Thus, the result holds.
\proofend 

We are now going to 
establish the following conclusion.

\prop{pro4-new}
{
	For any $m\in \N$ and 
	$t\in \{0,1,2,3\}$, 
	$s\in \A$ holds 
	for every integer $s$ 
	 satisfying 
	 \(m+2\le s\le3m+t-2\).
}

\myproof 
We will complete the proof in the four cases $t=0,1,2$ or $3$.
For each case, we 
apply Lemma~\ref{se4-l1}
to prove a sequence of auxiliary claims.

When $m=1$ and $t\in\{0,1\}$, there is no integer $s$
satisfying $m+2\le s\le 3m+t-2$, so the assertion holds
trivially. It remains to consider $m\ge2$, or $m=1$
and $t\in\{2,3\}$.

The following table summarizes the integer ranges supplied
by the claims below. 
All intervals in the table are integer intervals.

\vspace{2 mm}

\begin{tabularx}
{0.9\textwidth}
{@{}cp{15mm}Xp{45mm}@{}}
	\toprule
	$t$ & $m$ & 
	Ranges proved &
	Supporting claims
	 \\
	\midrule
	0 & $m\ge2$ & $[m+1,\,3m-2]$ &
	Claim 1.4 \\
	1 & $m\ge2$ & $[m,\,2m]\cup[2m,\,3m-1]$ &
		Claims 2.3 and 2.2 \\
	2 & $m\ge1$ & $[m,\,2m+2]\cup[2m+1,\,3m]$ & 
	Claims 3.2 and 3.3 \\
	3 & $m\ge1$ & $[m+1,\,2m+3]\cup[2m,\,3m+1]$ & 
		Claims 4.2 and 4.3 \\
	\bottomrule
\end{tabularx}

\vspace{2 mm}

In the remainder of this proof, 
any conclusion 
of the form $\Gamma\Bra{x_1^{[q_1]},
	x_2^{[q_2]}, x_3^{[q_3]}}\in \A$
is obtained 
by applying
the result of (\ref{se4-e0}),
and we omit further references to (\ref{se4-e0}). 
In every application 
below, the displayed multiplicities sum to at most
\(m+1\); the denominators and arguments are positive in the indicated range,
and the displayed value of \(\Gamma\) is nonnegative.

\incase $t=0$. Then $k=4m$.

There is no integer $s$
satisfying $m+2\le s\le 3m-2$
when $m=1$. 
Thus, we may assume that 
 $m\ge 2$ in 
 Case 1.

\noindent{\bf Claim 1.1}: 
$\frac{2(i+1)m}{i}\in \A$ for every $i\in \N$.

Let $a_1(i)\coloneq\frac{2(i+1)m}{i}$.
By Lemma~\ref{se4-l1}, $k=4m\in \A$. 
Thus, $a_1(1)\in \A$.
For any $i\in \N$, if 
$a_1(i)\in \A $,
then 
\equ{p0-cl1-1}
{
	a_1(i+1)=
\frac{2(i+2)m}{i+1}
=4m\Bra{1-\frac{m}{\frac{2(i+1)m}{i}}}
=4m\Bra{1-\frac{m}
	{a_1(i)}}
=\Gamma(a_1(i)^{[m]})\in \A.
}
Thus, $a_1(i)\in \A$ for 
every $i\in \N$
and 
Claim 1.1 holds.

\noindent{\bf Claim 1.2}: 
$2m\in \A$
and 
$2(m-1)\in \A$.

By Claim 1.1, 
$2(m+1)=a_1(m)\in \A$.
Thus,
\aln{
2m
&=4m\Bra{1-\frac{m+1}{2(m+1)}}
=\Gamma(a_1(m)^{[m+1]})
\in \A,\\
2(m-1)
&=4m\Bra{1-\frac{m+1}{2m}}
=\Gamma((2m)^{[m+1]})
\in \A.
}
 Claim 1.2 holds.

\noindent{\bf Claim 1.3}: 
$\frac{2m(m-i-1)}{m-i} \in \A$ for all integers 
$i$ with $0\le i\le m-2$.

Let $b_1(i)\coloneq\frac{2m(m-i-1)}{m-i}$
for $0\le i\le m-2$.
By Claim 1.2, 
$2(m-1)\in \A$.
Thus, $b_1(0)=2(m-1)\in \A$.
For any integer $i$
with $1\le i\le m-2$,  
if $b_1(i-1)\in \A$, then 
\equ{p0-cl3-1}
{
b_1(i)
=4m\Bra{1-\frac{m}{\frac{2m(m-i)}{m-i+1}}}
=4m\Bra{1
-\frac{m}{b_1(i-1)}}
=\Gamma(b_1(i-1)^{[m]})
\in \A.
}
Thus, Claim 1.3 holds.

\noindent{\bf Claim 1.4}: 
$s\in  \A$
for every integer $s$ 
with $m+1\le s\le 3m-2$.

For every integer $i$ with 
$0\le i\le m-2$,
by Claim 1.3, $b_1(i)=\frac{2m(m-i-1)}{m-i}\in \A$, implying that 
	$$
2m+i
=4m\Bra{1-\frac{m-i-1}
	{b_1(i)}-\frac i{4m}}
=\Gamma(b_1(i)^{[m-i-1]},k^{[i]})
\in \A.
$$
Hence $s\in \A$ for 
$2m\le s\le 3m-2$,
and in particular, 
$2m\in \A$.
For any integer $i$ with $0\le i\le m-2$, 
if $2m-i\in \A$, 
 then
 \aln{
 c_1(i)&\coloneq\frac{4m (m-i-1)}{2 m-i}
 =4m \Bra{1-\frac{m+1}{2m-i}}
 =\Gamma((2m-i)^{[m+1]})\in \A,
\\
 2 m-i-1
 &=4m\Bra{1
 -\frac{m-i-1}
 {c_1(i)}
 -\frac i{2m} -\frac 1{4m}}
 =\Gamma(c_1(i)^{[m-i-1]},(2m)^{[i]},k)\in \A.
}
 Thus,  $2m-i\in \A$
 for every integer $i$ with $0\le i\le m-1$.
 Claim 1.4 holds.

By Claim 1.4,
the conclusion holds in Case 1.

\incase $t=1$. Then $k=4m+1$.

There is no integer $s$
satisfying $m+2\le s\le 3m-1$
when $m=1$. 
Thus, we assume that 
$m\ge 2$ in Case 2.

\noindent{\bf Claim 2.1}: 
For every $i\in \N$,
$\frac{(i+1)(4m+1)}{2i+1} \in \A$.

Let $a_2(i)\coloneq\frac{(i+1)(4m+1)}{2i+1}$ and  $b_2(i)\coloneq\frac{(2i+1)m}
{i}$.
By Proposition~\ref{le8-1}, 
$b_2(1)=3m\in \A$.
For any $i\in \N$, if $b_2(i)\in\A$,
then
\aln{
a_2(i)&=\frac{(i+1)(4m+1)}{2i+1}
=(4m+1)\Bra{1-\frac{m}
	{b_2(i)} }
	=\Gamma(b_2(i)^{[m]})
	\in \A,
\\
b_2(i+1)
&=\frac{(2i+3)m}{(i+1)}
=(4m+1)\Bra{1-\frac{m}
	{a_2(i)}
	-\frac 1{4m+1}}
=\Gamma(a_2(i)^{[m]}, k)\in \A.
}
Thus, Claim 2.1 follows.

\noindent{\bf Claim 2.2}: 
$3 m-i\in \A$ holds
for every $i\in \brk{m}$.

By Claim 2.1, 
$a_2(i)=
\frac{(i+1)(4m+1)}{2i+1} \in \A$
for every $i\in \N$.
If $i\in \brk{m}$, then 
 \equ{p1-cl2-1}
{
	3 m-i
	=(4m+1)\Bra{1-\frac{i+1}
		{a_2(i)}
		-\frac{m-i}{4m+1}}
	=\Gamma(a_2(i)^{[i+1]}, k^{[m-i]})\in \A.
}
Claim 2.2 follows.

\noindent{\bf Claim 2.3}: 
For every $i\in \{0\}\cup \brk{m}$,
${2m-i}\in \A$.

Let $c_2(i)\coloneq{2m-i}$.
By Claim 2.2, 
$c_2(0)=2m\in \A$.
Note that $c_2(1)=2m-1$.
As $m\ge 2$ in Case 2, 
by definition, 
\aln{
d_2&\coloneq\frac{(m-1)(4m+1)}{2m}
=(4m+1)\Bra{1-\frac{m+1}{2m}}
=\Gamma((2m)^{[m+1]})\in \A,
\\
c_2(1)
&={2m-1}
=(4m+1)
\Bra{1-\frac{m-1}
	{d_2}
	 - \frac{2}{4m+1}}
=\Gamma(d_2^{[m-1]}, k^{[2]})\in \A.
}
Hence Claim 2.3 holds for 
$i\in \{0,1\}$.
It suffices to show that for
every integer $i$ with 
$0\le i\le m-2$, 
$c_2(i)\in \A$ always implies that 
$c_2(i+2)\in \A$.
If $c_2(i)\in \A$, then 
 \equ{p1-cl3-3}
{
e_2(i)\coloneq\frac{(m-i-1)(4m+1)}{2m-i}
=(4m+1)\Bra{1-\frac{m+1}{2m-i}}
=\Gamma(c_2(i)^{[m+1]})\in \A.
}
By Claim 2.1, 
$a_2(i+1)
=\frac{(i+2)(4m+1)}{2i+3} \in \A$.
Thus,
$$
2m-i-2=(4m+1)
\Bra{1-\frac{m-i-1}
	{e_2(i)}
-\frac{i+2}
{a_2(i+1)}
}
=\Gamma(e_2(i)^{[m-i-1]},a_2(i+1)^{[i+2]})\in \A,
$$
implying that $c_2(i+2)\in \A$.
Thus, Claim 2.3 holds.

By 
Claims 2.2 and 2.3, 
the conclusion holds in Case 2.

\noindent {\bf Case 3}: $t=2$.
Then $k=4m+2$.

In this case, we will 
show that $s\in \A$
for each integer with $m+2\le 
s\le 3m$.

\noindent{\bf Claim 3.1}: 
$\frac {(2i+1)m+i}{i}\in \A$
for every $i\in \N$.

Let $a_3(i)\coloneq\frac {(2i+1)m+i}{i}$
and 
$b_3(i)\coloneq\frac {(2i+1)m+i+1}{i}$. 
As
$a_3(1)
={3m+1}$
and 
$b_3(1)=3m+2$,
by Proposition~\ref{le8-1},
both $a_3(1)$ and $b_3(1)$ 
belong to $\A$.
Assume inductively that 
both $a_3(i)$ and $b_3(i)$
belong to $\A$ for some $i\in \N$.
Then 
\aln
{
c_3(i)&\coloneq\frac{(4 m+2)m (i+1)}{i(2m+1)+m}
=(4m+2)\Bra{1-\frac{m+1}
	{a_3(i)}}
=\Gamma(a_3(i)^{[m+1]})\in \A,
\\
b_3(i+1)&=
\frac{(2i+3)m+i+2}{(i+1)}
=(4m+2)\Bra{1-\frac{m}
	{c_3(i)}}
=\Gamma(c_3(i)^{[m]})\in \A,
\\
	d_3(i)&\coloneq
	\frac{(4 m+2) (i+1)(m+1)}{i(2m+1) +m+1}
	=(4m+2)
	\Bra{1-\frac{m}{b_3(i)
	}}
	=\Gamma(b_3(i)^{[m]})\in \A,
\\
a_3(i+1)&=
\frac{(2i+3)m+(i+1)}
{(i+1)}
=(4m+2)\Bra{1
	-\frac{m+1}{d_3(i)
	}
}
=\Gamma(d_3(i)^{[m+1]})\in \A.
}
Thus, 
$a_3(i)\in \A$ and $b_3(i)\in \A$
for every $i\in \N$, and 
Claim 3.1 holds.

\noindent{\bf Claim 3.2}: 
$2m-i
\in \A$ holds for every integer $i$ with 
$-2\le i\le m$.

Let $e_3(i)\coloneq{2m-i}$.
Then 
$e_3(-2)=2m+2=a_3(m)
 \in \A$
by Claim 3.1.
By definition,
\equ{p2-cl2-1}
{
e_3(-1)=2m+1
=(4m+2)\Bra{1-\frac{m+1}{2m+2}}
=\Gamma(e_3(-2)^{[m+1]})\in \A.
}
Thus, $e_3(i)\in \A$
for each $i\in \{-2,-1\}$.
For any integer $i$  
with $-2\le i\le m-2$,
if $e_3(i)\in  \A$,
then
\aln{
f_3(i)&\coloneq \frac{(m-i-1)
	(4 m+2)}{2 m-i}
=(4m+2)\Bra{1-\frac{m+1}{2m-i}}
=\Gamma(e_3(i)^{[m+1]})
\in \A, \\
 	2m-i-2
 	&=(4m+2)\Bra{1
 		-\frac{m-i-1}{f_3(i)
 			}
 			-\frac{i+2}{2m+1}
 	}
 =\Gamma(f_3(i)^{[m-i-1]}, e_3(-1)^{[i+2]})\in \A,
}
 implying that $e_3(i+2)\in \A$.
 Thus, Claim 3.2 holds.
 
\noindent{\bf Claim 3.3}:  
 ${s}\in \A
 $ for every integer $s$ with 
 $2m+1\le s\le 3m$.

By Claim 3.2, $2m\in \A$
and $2m+1\in \A$.
For every integer $i$ with 
$0\le i\le m-2$, 
by  
the proof of Claim 3.2, 
$f_3(i)=
\frac{ (m-i-1)(4 m+2)}{2 m-i}\in 
\A$, implying that 
 \equ{p2-cl3-1}
{
2 m+2+i
=(4m+2)
\Bra{1-
\frac{m-i-1}
{f_3(i)
} 
}
=\Gamma(f_3(i)^{[m-i-1]})\in \A.
}
Hence Claim 3.3 holds.

By  
Claims 3.2 and 3.3, the conclusion holds in Case 3.

\noindent {\bf Case 4}: $t=3$.
Then $k=4m+3$.

In this case, we will 
show that $s\in \A$
for each integer with $m+2\le 
s\le 3m+1$.

\noindent{\bf Claim 4.1}: 
$\frac{(i + 1)(4m + 3)}{2i + 1}\in \A$ for $0\le i\le m+1$
and $\frac{(2i+1)(m+1)}{i}
\in \A$ for $i\in  \brk{m+1}$.

Let $a_4(i)\coloneq\frac{(i + 1)
	(4m + 3)}{2i + 1}$
and $b_4(i)\coloneq\frac{(2i+1)(m+1)}
{i}$.
Note that $a_4(0)=4m+3\in \A$
by Lemma~\ref{se4-l1}, and 
\aln{
b_4(1)&=3(m+1)
=(4m+3)\Bra{1-\frac{m}{4m+3}}
=\Gamma(k^{[m]})\in \A,
\\
a_4(1)&=\frac{2(4 m+3)}{3}
=(4m+3)\Bra{1-\frac{m+1}{3m+3}}
=\Gamma(b_4(1)^{[m+1]})\in \A.
}
Assume
inductively that 
$a_4(i-1), a_4(i),b_4(i)\in \A$
for some $i\in \brk{m}$.  
Let 
$	c_4(i):=
\frac{2 m i+i+m}{i}$. Then
\aln{
	c_4(i)&=
(4m+3)\Bra{1
	-\frac{i+1}
	{a_4(i)}
-\frac{m-i}{a_4(i-1)}
}
=\Gamma(a_4(i)^{[i+1]}, a_4(i-1)^{[m-i]})\in \A,
\\
	d_4(i)&\coloneq
\frac{(i+1)m(4 m+3)}{2 m i+i+m}
	=(4m+3)\Bra{1
	-\frac{m+1}{c_4(i)}
}
	=\Gamma(c_4(i)^{[m+1]})
	\in \A,
\\
b_4(i+1)
&=\frac{(2i+3) (m+1)}
{(i+1)}
=(4m+3)\Bra{1-\frac{m}
{d_4(i)}
}
=\Gamma(d_4(i)^{[m]})\in \A,
\\
a_4(i+1)
&=\frac{ (i+2)(4 m+3)}{2 i+3}
=(4m+3)
\Bra{1-\frac{m+1}
	{b_4(i+1)}
}
=\Gamma(b_4(i+1)^{[m+1]})
\in \A.
}
Thus, Claim 4.1 holds.

\noindent{\bf Claim 4.2}: 
${2m-i}\in \A$ holds 
for every integer $i$ with 
$-3\le i\le m-1$.

Let $e_4(i)\coloneq 2m-i$.
By Claim 4.1, $e_4(-3)=2m+3=b_4(m+1)\in \A$.
By definition,
 \aln
{
	 f_4&\coloneq
	 2m+2=(4m+3)
	\Bra{1-\frac{m+1}
	{a_4(m)}
	}
	=\Gamma(a_4(m)^{[m+1]})
	\in \A,
\\
 g_4&\coloneq	\frac{4m+3}{2}
 	=(4m+3)\Bra{1-\frac{m+1}{2m+2}}
 	=\Gamma(f_4^{[m+1]})\in \A,
\\
 	h_4&\coloneq 2m+1
 	=(4m+3)
 	\Bra{1-\frac{m+1}
 		{g_4}}
 	=\Gamma(g_4^{[m+1]})\in \A,
\\
 	p_4&\coloneq \frac{m(4m+3)}{2m+1}
 	=(4m+3)
 	\Bra{1-\frac{m+1}{2m+1}}
 	=\Gamma(h_4^{[m+1]})\in \A,
\\
 	q_4&\coloneq 2m
 	=(4m+3)
 	\Bra{1-\frac{m}
 		{p_4}
 		-\frac 1{g_4}}
 	=\Gamma(p_4^{[m]}, g_4)\in \A,
}
  implying that $e_4(i)\in\A$
  for each $i\in \{-2,-1,0\}$.
  Thus, Claim 4.2 holds when $m=1$.
  
  Now assume that $m\ge 2$.  
  It remains to show that for each
  integer $i$ with $0\le i\le m-2$, 
  $e_4(i)\in \A$ always
  implies that 
   $e_4(i+1)\in \A$.
Assume that $e_4(i)\in \A$. Then, 
  \equ{p3-cl2-6}
{
u_4(i)\coloneq
\frac{ (m-i-1)(4 m+3)}{2 m-i}
=(4m+3)
\Bra{1-\frac{m+1}{2m-i}}
=\Gamma(e_4(i)^{[m+1]})
\in \A.
}
As $g_4=\frac{4 m+3}{2}
\in \A$, 
$$
e_4(i+1)=2 m-i-1
=(4m+3)
\Bra{1-\frac{m-i-1}{u_4(i)}
	-\frac{i+2}{\frac{4 m+3}{2}}}
=\Gamma(u_4(i)^{[m-i-1]},g_4^{[i+2]})\in \A.
$$
Hence Claim 4.2 follows.

\noindent{\bf Claim 4.3}: 
$s\in \A$ holds for every 
integer 
 $s$ with 
$2m\le s\le 3m+1$.

Let $v_4(i)\coloneq2m+i$.
By Claim 4.2, $v_4(0),v_4(1)\in  \A$.
For $0\le i\le m-1$,  
by Claim 4.1, 
$a_4(m-i)=\frac{(m-i+1)
	(4 m+3)}{2(m-i)+1}\in\A$,
and thus,
  \equ{p3-cl3-1}
{
2 m+2+i
=(4m+3)
\Bra{1-\frac{m+1-i}
{a_4(m-i)}
		-\frac i{4m+3}}
		=\Gamma(a_4(m-i)^{[m+1-i]},k^{[i]})\in \A,
}
implying that $v_4(i+2)\in \A$.
Hence
$s\in \A$
for each integer $s$ with 
$2m\le s\le 3m+1$, and  Claim 4.3 holds.

By 
Claims 4.2 and 4.3, 
the result holds in Case 4.

Therefore, Proposition~\ref{pro4-new} holds.
\proofend 

\vspace{0.3 cm}
\noindent {\it Proof} of Theorem~\ref{main-th2}.
The necessity follows from 
Theorem~\ref{known0}.
Now we prove the sufficiency.

Let $k=4m+t$ with $k\ge 2$, 
where $m, t\ge 0$ are integers
and $0\le t\le 3$, 
and let $r_0=\Floor{\frac 14 k} +2=m+2$.
For every integer \(r\), the inequality \(\frac{k}{4}+1<r\) is equivalent to \(r\ge r_0\).

If $m=0$, then $k\in\{2,3\}$ and Proposition~\ref{le8-1}(i) gives $s\in\mathcal A$ for every integer $s$ with $m+2\le s\le k$. If $m\ge1$, Proposition~\ref{pro4-new} supplies the integers satisfying $m+2\le s\le3m+t-2$, while Proposition~\ref{le8-1}(ii) supplies those satisfying $3m+t-1\le s\le4m+t=k$. Together these cover every integer $s$ with $m+2\le s\le k$. Consequently,
\[
\frac{s}{\sqrt{k}}\in\mathscr W_{r_0}(\sqrt{k})
\subseteq\mathscr W_s(\sqrt{k})
\qquad(r_0=m+2).
\]

For every integer $s$ in $[r_0,k]$,  since 
$
s\times 
\left (\frac{s}{\sqrt{k}}\right )^{-1}
=\sqrt{k}, 
$
Theorem~\ref{main-th01}
implies the existence of 
a tree $T$ with $\Delta(T)=s$ and $\rho(T)=\sqrt k$.
\proofend

\section{Further results
	and questions
\label{problem} 
}

Recall that $\Rho$ is the set of 
values $\rho(T)$ attained by trees $T$ 
of order at least $2$.
Problem~\ref{conj1} asks 
whether a tree $T$ exists
with $\rho(T)=\alpha$ 
and $\Delta(T)=r$ 
for any given $\alpha\in \Rho\setminus \{1\}$ and integer 
$r$ satisfying 
$\frac 14 \alpha^2+1<r\le \alpha^2$.
Theorem~\ref{main-th2}
gives a positive answer to  Problem~\ref{conj1}
whenever 
$\alpha=\sqrt k$ for 
an integer $k\ge 2$.

For any integer $q\ge2$ and  
$k_1,\ldots,k_q\in \N$, the eigenvalue formula for Cartesian
products (see \cite[Sections~1.4.2 and~1.4.6]{Brouwer})
gives
\[
\rho\bigl(K_{1,k_1}\square\cdots\square K_{1,k_q}\bigr)
=\sum_{i=1}^{q}\sqrt{k_i}.
\]

Thus the construction realizes the sum
as the spectral radius of a finite graph containing a $4$-cycle
and leads to
the question of whether it can also be realized as
the spectral radius of a tree.

Since $\rho(K_{1,k})=\sqrt{k}$, we have
$\sqrt{k}\in\mathbb P$ for every $k\in\mathbb N^+$.
The following result shows that 
the sum of the square roots of any two numbers in 
$\N$ also belongs to $\Rho$.

\prop{prop5-1}
{
For any $k_1,k_2\in\N$, the number $\sqrt{k_1}+\sqrt{k_2}$ belongs to $\Rho$.
If $k_1\ne k_2$, it is realized by a tree of maximum degree $k_1+k_2$.
}

\proof 
Assume that $k_1\le k_2$. 
If $k_1=k_2$, then 
$\sqrt{k_1}+\sqrt{k_2}=
2\sqrt{k_1}=\rho(K_{1,4k_1})\in 
\Rho$.
If $k_1<k_2$, put 
$\alpha=\sqrt{k_1}+\sqrt{k_2}$ and $r=k_1+k_2$.
Then, 
\aln
{
c_1:&=\alpha
-\frac{k_2-k_1}
{\sqrt{k_1}+\sqrt{k_2}}
=2\sqrt {k_1},\\
c_2:&=\alpha
-\frac{2k_1}{c_1}
=\sqrt {k_2},\\
c_3:&=\alpha-\frac{k_2}{c_2}
=\sqrt{k_1}.
}
All three outputs $c_1,c_2$ 
and $c_3$ are positive and the corresponding arities are less than $r$, so they belong to
$\scre_r(\alpha)$.
Since 
$$
\frac{k_1}{c_3}+\frac{k_2}{c_2}
=\sqrt {k_1}+\sqrt{k_2}
=\alpha, 
$$
Theorem~\ref{main-th01} 
gives a tree $T$ 
with $\Delta(T)=r=k_1+k_2$ and 
$\rho(T)=
\sqrt{k_1}+\sqrt{k_2}\in \Rho$.
\proofend 

This leads to the following question.

\begin{problem}\label{se5-p2}
For an integer \(q\ge3\), characterize the tuples
\((k_1,\ldots,k_q)\in(\mathbb N^+)^q\) for which
\[
\sum_{i=1}^{q}
\sqrt{k_i}\in\Rho.
\]
\end{problem}

Write each \(k_i=a_i^2d_i\), where \(a_i,d_i\in \N\)  and \(d_i\) is squarefree. 
If all
$d_i$ are equal to $d$, then
\[
\sum_{i=1}^{q}\sqrt{k_i}=A\sqrt d=\sqrt{A^2d}\in \Rho,
\qquad A=\sum_{i=1}^{q}a_i.
\]
If exactly two distinct squarefree parts $d,e$ occur, put
$A=\sum_{d_i=d}a_i$ and $B=\sum_{d_i=e}a_i$. Then $A,B\in\N$ and
\[
\sum_{i=1}^{q}\sqrt{k_i}=\sqrt{A^2d}+\sqrt{B^2e}\in \Rho
\]
by  Proposition~\ref{prop5-1}. Thus the unresolved cases in Problem~\ref{se5-p2} have
at least three distinct squarefree parts.

We next consider averages of square roots of positive integers. The following
result characterizes when the average of two such square roots belongs to
$\Rho$.

\prop{prop5-2}
{
	For any $p,q\in\N$, the number $\frac 12 \Bra{\sqrt{p}+\sqrt{q}}$ belongs to $\Rho$
	if and only if 
	$p\equiv q\pmod4$.
}

\proof
{\bf Sufficiency}. 
Assume that $p\ge q$
and $p-q=4D$. 
If $p=q$,  
then $\alpha=\sqrt q=\rho(K_{1,q})$, 
implying that $\alpha\in \Rho$.

Now assume that $p>q$. 
Put $r=\max\{q,D\}+1$ and 
$\alpha=
\frac{\sqrt p+\sqrt{q}}2$.
Then 
$D>0$ and 
\aln{
c_1:&=\alpha
-\frac D{\alpha}
=\frac{\sqrt p+\sqrt q}2
-\frac{\sqrt p-\sqrt q}2
=\sqrt{q},\\
c_2:&=\alpha
-\frac{q}{c_1}
=\frac {\sqrt p+\sqrt q}2-
\frac{q}{\sqrt q}
=\frac{\sqrt p-\sqrt{q}}2.
}
Since both $c_1$ 
and $c_2$ are positive and the corresponding arities are less than $r$,  
by the definition of $\scre_r(\alpha)$, 
both belong to
$\scre_r(\alpha)$.
Since $p> q$, $c_2>0$ and 
$$
0=\frac{\sqrt p+\sqrt{q}}2
-\frac{\sqrt p+\sqrt{q}}2
=\alpha-\frac {\frac{p-q}{4}}{\frac{\sqrt p-\sqrt{q}}2}
=\alpha-\frac D
{c_2}
\in \scre_r(\alpha).
$$
Hence $0\in \scre_r(\alpha)$.
By Theorem~\ref{main-th0}, 
there exists a tree $T$ with $\Delta(T)\le r$ and $\rho(T)=\alpha$,
implying that  $\alpha=\frac{\sqrt p+\sqrt q}2\in \Rho$.
The sufficiency holds.

\noindent 
\textbf{Necessity.} 
Assume that $\alpha=\frac{\sqrt p+\sqrt q}2\in \Rho$.
The eigenvalues of an integer matrix are algebraic
integers. Hence $\alpha$ is an algebraic
integer. Write $p=a^2d$ and $q=b^2e$, where $a,b\in\N$ and $d,e$ are
positive squarefree integers.

\begin{enumerate}
[itemsep=-1mm, label=(\roman*)]
	\item 
If $d=e=1$, then $\alpha=(a+b)/2$ is a rational algebraic integer, and
therefore an integer. Thus $a,b$ have the same parity, giving
$p\equiv q\pmod4$.

	\item 
If $d=e>1$, the monic minimal polynomial of $\alpha$ over
$\Q$
is
\[
x^2-\frac{d(a+b)^2}{4}.
\]
Its coefficients must be integers. Since $d$ is squarefree, 
$4\not \mid d$.
Thus, 
$d(a+b)^2/4\in\N$ forces $a+b$ to be even. Again $a,b$ have the same
parity, and $p\equiv q\pmod4$.

	\item 
If $d\ne e$ and one of $d,e$ equals $1$, assume $d=1$.
The monic minimal polynomial is
\[
x^2-a x+\frac{p-q}{4}.
\]
Its constant coefficient is an integer, so $4\mid(p-q)$.

	\item 
Finally, if $d,e>1$ and $d\ne e$, the field
$\mathbb Q(\sqrt d,\sqrt e)$ has degree four. The four conjugates of
$\alpha$ are 
$(\pm\sqrt p\pm\sqrt q)/2$, with independent signs, and are
distinct. The monic minimal polynomial is therefore
\[
x^4-\frac{p+q}{2}x^2
+\frac{(p-q)^2}{16}.
\]
Integrality of its constant coefficient gives $16\mid(p-q)^2$, hence
$4\mid(p-q)$. 
\end{enumerate}
These cases prove necessity.
\proofend

This leads to the following question.

\begin{problem}\label{se5-p3}
	For an integer \(q\ge3\), characterize the tuples
	\((k_1,\ldots,k_q)\in(\mathbb N^+)^q\) for which
	\[
	\frac 1q\sum_{i=1}^{q}
	\sqrt{k_i}\in\Rho.
	\]
\end{problem}

A parallel recursive characterization for Laplacian spectral radii of trees
was developed in 
\cite{Zhang:2025a}.

\vspace{3mm} 

\noindent {\bf Conflict of interest}: The authors declare that they have no known competing financial interests or personal
relationships that could have appeared to influence the work reported in this paper.

\noindent {\bf Data availability statement}:  Not applicable.

\section*{Acknowledgments}
This research is supported by 
NSFC (Nos. 12101347 and 12371340), 
NSF of Shandong Province (No. ZR2021QA085)  and 
the Ministry of Education,
Singapore, under its Academic Research Tier 1 (RG19/22). 
Any opinions,
findings and conclusions or recommendations expressed in this
material are those of the authors and do not reflect the views of the
Ministry of Education, Singapore.

\def \LAA {{\it Linear Algebra Appl.} }
\def \JCTB {{\it J. Combin. Theory, Ser. B} }
\def \JGT {{\it J. Graph Theory} }

\end{document}